\documentclass[11pt]{article}
\usepackage{cite}
\usepackage{mathrsfs}
\usepackage{amsfonts}
\usepackage{amsmath}
\usepackage{amsfonts,amssymb,color}
\usepackage{dsfont}
\usepackage{curves}
\usepackage{mathrsfs}
\usepackage{pifont}
\usepackage{amssymb}
\usepackage{latexsym,amsmath,amssymb,amsfonts,epsfig,graphicx,cite,psfrag}
\usepackage{eepic,color,colordvi,amscd}

\newtheorem{theorem}{Theorem}[section]

\newtheorem{lemma}{Lemma}[section]

\newtheorem{problem}{Problem}[section]
\newtheorem{claim}{Claim}[section]

\newcommand{\qed}{\hfill\rule{0.5em}{0.809em}}

\def\emptyset{\mbox{{\rm \O}}}

\def\qed{\hfill \rule{4pt}{7pt}}

\def\pf{\noindent {\it Proof. }}

\begin{document}
	
	\title{The optimal binding function for (cap, even hole)-free graphs}
	\author{Ran Chen$^{1,}$\footnote{Email: 1918549795@qq.com},  \;Baogang  Xu$^{1,}$\footnote{Email: baogxu@njnu.edu.cn. Supported by 2024YFA1013902 and NSFC 11931006}, \;\; Yian Xu$^{2,}$\footnote{Corresponding author. Email: yian$\_$xu@seu.edu.cn. Supported by NSFC 12471322}\\\\
		\small $^1$Institute of Mathematics, School of Mathematical Sciences\\
		\small Nanjing Normal University, 1 Wenyuan Road,  Nanjing, 210023,  China\\
		\small $^2$School of Mathematics, Southeast University, 2 SEU Road, Nanjing, 211189, China}
\date{}
	
	\maketitle
	
\begin{abstract}
A {\em hole} is an induced cycle of length at least 4, an {\em even hole} is a hole of even length, and a {\em cap} is a graph obtained from a hole by adding an additional vertex which is adjacent exactly to two adjacent vertices of the hole. A graph $G$ obtained from a graph $H$ by blowing up all the vertices into cliques is said to be a clique blowup of $H$. Let $p, q$ be two positive integers with $p>2q$, let $F$ be a triangle-free graph, and let $G'$ be a clique blowup  of $F$ with $\omega(G')\leq\max\{\frac{2q(p-q-2)}{p-2q}, 2q\}$. In this paper, we prove that for any clique blowup $G$ of $F$, $\chi(G)\leq\lceil\frac{p}{2q}\omega(G)\rceil$ if and only if $\chi(G')\leq\lceil\frac{p}{2q}\omega(G')\rceil$. As its consequences, we show that every (cap, even hole)-free graph $G$ satisfies $\chi(G)\leq\lceil\frac{5}{4}\omega(G)\rceil$, which affirmatively answers a question of Cameron {\em et al.} \cite{CdHV2018}, we also show that every (cap, even hole, 5-hole)-free graph $G$ satisfies  $\chi(G)\leq\lceil\frac{7}{6}\omega(G)\rceil$, and the bound is reachable.
\begin{flushleft}
\end{flushleft}
		
\end{abstract}
	
\newpage
	
\section{Introduction}

All graphs in this paper are finite and simple. We follow \cite{BM08} for undefined notations and terminologies. Let $G=(V, E)$ be a graph. Let $N_G(v)$ be the set of neighbors of $v$, $d_G(v)=|N_G(v)|$. If it does not cause any confusion, we usually omit the subscript $G$. Let $X$ be a subset of $V(G)$. We use $G[X]$ to denote the subgraph of $G$ induced by $X$.

We say that a graph $G$ contains a graph $H$ if $H$ is isomorphic to an induced subgraph of $G$, and say that $G$ is $H$-{\em free} if it does not contain $H$.
For a family $\{H_1,H_2,\cdots\}$ of graphs, $G$ is $(H_1, H_2,\cdots)$-free if $G$ is
$H$-free for every $H\in \{H_1,H_2,\cdots\}$.

A {\em clique} (resp. {\em stable set}) of $G$ is a set of mutually adjacent
(resp. non-adjacent) vertices in $G$. The {\em clique number} (resp. {\em stability number}) of $G$, denoted by $\omega(G)$ (resp. $\alpha(G)$), is the maximum size of a clique (resp. stable set) in $G$.

Let $k$ be a positive integer. A $k$-{\em coloring } of $G$ is a function $\phi: V(G)\rightarrow \{1,\cdots,k\}$ such that $\phi(u)\ne \phi(v)$ if $uv\in E(G)$. The {\em chromatic number} $\chi(G)$ of $G$ is the minimum number $k$ for which $G$ has a $k$-coloring. A graph is {\em perfect} if all its induced subgraphs $H$ satisfy $\chi(H)=\omega(H)$. A {\em hole} is an induced cycle of length at least 4, an {\em even hole} (resp. {\em odd hole}) is a hole of even (resp. odd) length. In 2006, Chudnovsky {\em et al} \cite{CRST2006} proved the {\em Strong Perfect Graph Theorem} stating that a graph $G$ is perfect if and only if $G$ does not contain an odd hole or its complement. However, for any two positive integers $k$ and $\ell$, Erd\H{o}s \cite{E1959} showed that there exists a graph $G$ with $\chi(G)\geq k$ whose shortest cycle has length at least $\ell$, and so $\chi(G)-\omega(G)$ may be arbitrarily large.

In 1975,  Gy\'{a}rf\'{a}s \cite{G75} proposed the following important concept. Let ${\cal G}$ be a family of graphs.  If there exists a function $f$ such that $\chi(G)\leq f(\omega(G))$ for all graphs $G$ in ${\cal G}$, then we say that ${\cal G}$ is $\chi$-{\em bounded}, and call $f$ a {\em binding function} of ${\cal G}$.

Dirac \cite{D61} showed that every chordal graph (i.e., a hole-free graph) has a simplicial vertex (i.e., whose neighbor set is a clique). This implies that $\chi(G)=\omega(G)$ if $G$ is a hole-free graph. Vizing \cite{V65} showed that each line graph $G$ satisfies  $\chi(G)\leq\omega(G)+1$. In 2016, Scott and Seymour \cite{SS2016} showed that $\chi(G)\leq\frac{2^{2^{\omega(G)+1}}}{48(\omega(G)+1)}$ if $G$ is odd hole-free. Then, Chudnovsky, Scott and Seymour \cite{CSS2017} proved that (hole of length at least $\ell$)-free graphs are $\chi$-bounded, and later Chudnovsky, Scott, Seymour and Spirkl \cite{CSSS2020} proved that (odd hole of length at least $\ell$)-free graphs are $\chi$-bounded. These three results confirm three long-standing conjectures of Gy\'{a}rf\'{a}s \cite{G87}.

The class of even hole-free graphs has been studied extensively, and interest readers are referred to \cite{SR2019} and \cite{V2010} for comprehensive surveys of even hole-free graphs. In 2001, Reed \cite{R2001} conjectured that every even hole-free graph has a bisimplicial vertex (i.e., whose neighbor set is the union of two cliques). In 2023, Chudnovsky and Seymour \cite{CS2023} confirmed this conjecture, and this implies that every even hole-free graph $G$ satisfies $\chi(G)\leq2\omega(G)-1$. This is a nice chromatic bound compared with the bound for odd hole-free graphs given by Scott and Seymour \cite{SS2016}. Cameron {\em et al.} \cite{CCH18} proved that every (pan, even hole)-free graph $G$ satisfies  $\chi(G)\leq\frac{3}{2}\omega(G)$, where a {\em pan} is a graph consisting of a hole and an additional vertex adjacent to exactly one vertex of the hole. Kloks {\em et al.} \cite{KMV2009} proved that $\chi(G)\leq \omega(G)+1$ if $G$ is (diamond, even hole)-free, where {\em diamond} is a graph consisting of two triangles sharing exactly one edge, and then Fraser {\em et al.} \cite{FHH18} extended the bound to (kite, even hole)-free graphs, where a {\em kite} is a graph consisting of a diamond and a new vertex adjacent to exactly one vertex of degree two of the diamond.

\begin{figure}[htbp]
	\begin{center}
		\includegraphics[width=10cm]{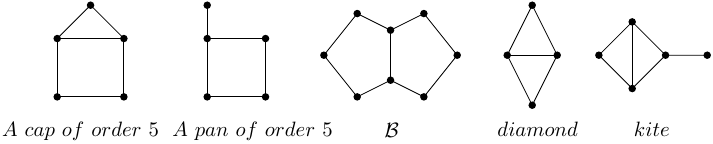}
	\end{center}
	\vskip -25pt
	\caption{Illustration of some forbidden graphs.}
	\label{fig-10}
\end{figure}

A {\em cap} is a graph consisting of a hole and an additional vertex which has exactly two adjacent neighbors in the hole. The class of (cap, odd hole)-free graphs has been studied extensively (see \cite{CCKV1999}). Since the complement of an odd hole of length at least 7 contains a cap of order 5, it is easy to see, from the Strong Perfect Graph Theorem, that (cap, odd hole)-free graphs are perfect. But the structure of (cap, even hole)-free graphs is very different from that of perfect graphs.

In 2018, Cameron, Silva, Huang and Vu\v{s}kovi\'c \cite{CdHV2018} showed that the minimum degree of an (cap, even hole)-free graph $G$ is at most $\frac{3}{2}\omega(G)-1$ and this implies that each (cap, even hole)-free graph $G$ satisfies that $\chi(G)\leq\frac{3}{2}\omega(G)$.  They also proposed the following problem.

\begin{problem}\label{pro}{\em \cite{CCH18}}
	 Is it true that $\chi(G)\leq\lceil\frac{5}{4}\omega(G)\rceil$ for every (cap, even hole)-free graph $G$?
\end{problem}

{\em Substituting} a vertex $v$ of a graph $G$ by a graph $H$ is an operation which creates a new graph with vertex set $V(H)\cup V(G-v)$ and edge set $E(G-v)\cup \{xy~|~x\in H, y\in N_{G}(v)\}$. When $H$ is a clique (not necessary nonempty), the substitution is said to be {\em blowing up} $v$ of $G$ into a clique. A graph $G$ obtained from a graph $H$ by blowing up all the vertices into cliques is said to be a {\em clique blowup } of $H$. Moreover, if each such a clique is nonempty, then it is said to be a {\em nonempty clique blowup} of $H$; if each such a clique has size $t$, then it is said to be a {\em $t$-clique blowup} of $H$ and denoted by $H^t$.

Let ${\cal B}$ be the graph consisting of two 5-holes sharing exactly one edge. Wu and Xu \cite{WX2019} proved that every (cap, even hole)-free graph has $\chi(G)\leq\lceil\frac{4}{3}\omega(G)\rceil$, and they showed that if $G$ is further ${\cal B}$-free, then $\chi(G)\leq\lceil\frac{5}{4}\omega(G)\rceil$. Later, Xu \cite{X2021} improved the bound $\lceil\frac{4}{3}\omega(G)\rceil$ to $\lceil\frac{9}{7}\omega(G)\rceil+1$. In this paper, we prove that every (cap, even hole)-free graph $G$ satisfies  $\chi(G)\leq\lceil\frac{5}{4}\omega(G)\rceil$. This affirmatively answers Problem~\ref{pro}. Moreover, we show that every (cap, even hole, 5-hole)-free graph $G$ satisfies  $\chi(G)\leq\lceil\frac{7}{6}\omega(G)\rceil$. Let $k$ be a positive integer, and let $G_1$ be a 5-hole and $G_2$ be a 7-hole. We have that $\chi(G_1^k)\geq\lceil\frac{|V(G_1^k)|}{\alpha(G_1^k)}\rceil=\lceil\frac{5k}{2}\rceil=\lceil\frac{5}{4}\omega(G_1^k)\rceil$ and $\chi(G_2^k)\geq\lceil\frac{|V(G_2^k)|}{\alpha(G_2^k)}\rceil=\lceil\frac{7k}{3}\rceil=\lceil\frac{7}{6}\omega(G_2^k)\rceil$. Therefore, both the bounds are reachable.

In the following Theorem~\ref{blowup}, we transfer Problem~\ref{pro} to graphs with bounded clique number. Then, we affirmatively answer Problem~\ref{pro}, in Theorem~\ref{main-1}, by first confirming it on graphs with small clique number.

\begin{theorem}\label{blowup}
	Let $p, q$ be two positive integers with $p>2q$, let $F$ be a triangle-free graph, and let $G'$ be a clique blowup  of $F$ with $\omega(G')\leq\max\{\frac{2q(p-q-2)}{p-2q}, 2q\}$. Then for any clique blowup $G$ of $F$, $\chi(G)\leq\lceil\frac{p}{2q}\omega(G)\rceil$ if and only if $\chi(G')\leq\lceil\frac{p}{2q}\omega(G')\rceil$.
\end{theorem}

\begin{theorem}\label{main-1}
	Let $G$ be an (cap, even hole)-free graph. Then $\chi(G)\leq\lceil\frac{5}{4}\omega(G)\rceil$. Moreover, if $G$ is further $5$-hole-free, then $\chi(G)\leq\lceil\frac{7}{6}\omega(G)\rceil$. Both bounds are reachable.
\end{theorem}

In the following Section 2, we present some useful notations and lemmas. Then, we prove Theorems~\ref{blowup} and \ref{main-1} in Sections 3 and 4, respectively.

\section{Notations and Preliminary Results}

In this section, we introduce some notations and conclusions which are useful to our proofs.

Let $G$ be a graph, and let $f$ be an assignment of integers 0 and 1 to its edges. A subgraph $H$ of $G$ is said to be {\em odd} (resp. {\em even}) if $\sum_{e\in E(H)}f(e)$ is odd (resp. even). The graph $G$ is said to be {\em odd-signable}  if it has an assignment such that every induced cycle is odd, and is {\em even-signable} if it has an assignment such that every triangle is odd and every hole is even. It is certain that every even hole-free graph $G$ is odd-signable by assigning each edge of $G$ with integer 1.

Let $H$ be a hole of $G$, and let $xyz$ be a segment of $H$. An induced $xz$-path $P$ is called an {\em ear} of the hole $H$ if the internal vertices of $P$ are all in $V(G)\setminus V(H)$, and $G[(V(H)\setminus \{y\})\cup V(P)]$ is also a hole. We call $x$ and $z$ the {\em attachments} of the ear $P$ in $H$.

A {\em wheel} $(H,v)$ of $G$ consists of a hole $H$ together with a vertex $v$, which is called the {\em center} of the wheel, such that $v$ has at least three neighbors on $H$.


A graph $G$ is said to be obtained from a graph $F$ by an {\em ear addition} if (1) the vertices of $V(G)\setminus V(F)$ are those internal vertices
of an ear, say $P$ with attachments $x$ and $z$, of some hole $H$ in $F$, and (2) $G$ contains no edge between vertices of $V(P)\setminus \{x,z\}$ and $V(G)\setminus\{x,y,z\}$, where $x,y,z$ are three consecutive vertices of $H$. An ear addition is said to be {\em good} if

\begin{itemize}
	\item $y$ has an odd number of neighbors on $V(P)$, and
	\item $F$ contains no wheel $(H_1,v)$, where $x,y,z\in V(H_1)$ and $vy\in E(G)$, and
	\item $F$ contains no wheel $(H_2,y)$, where $x,z$ are both neighbors of $y$ in $H_2$.
\end{itemize}

The graph obtained from the complete bipartite graph $K_{4,4}$ by removing a perfect matching is called a {\em cube}.  A {\em clique cutset} of $G$ is a clique $K$ in $G$ such that $G-K$ has more components than $G$. A clique $K$ is called a {\em universal clique} of $G$ if $K$ is complete to $V(G)\setminus K$.

\begin{lemma}\label{odd-signable}{\em \cite{CCKV2000}}
Let $G$ be a triangle-free graph with at least three vertices that is not the cube and has no clique cutset. Then, $G$ is odd-signable if and only if it can be obtained, starting from a hole, by a sequence of good ear additions.
\end{lemma}

\begin{lemma}\label{cap}{\em \cite{CdHV2018}}
Let $G$ be a (cap, $4$-hole)-free graph that contains a hole and has no clique cutset. Let $F$ be any maximal triangle-free induced subgraph of $G$ with at least $3$ vertices and has no clique cutset. Then $G$ is obtained from $F$ by first blowing up vertices of $F$ into nonempty cliques, and then adding a universal clique.
\end{lemma}

\section{Proof of Theorem~\ref{blowup}}

Obviously, the necessity of Theorem~\ref{blowup} is correct, and so it is left to prove the sufficiency. Let $p,q$ be two positive integers such that $p>2q$, and let $F$ be a triangle-free graph. We prove by induction on $|V(G)|$ that $\chi(G)\leq\lceil\frac{p}{2q}\omega(G)\rceil$.
Let $G$ be a clique blowup of $F$ with $\omega(G)\geq\max\{\frac{2q(p-q-2)}{p-2q}, 2q\}+1$, and suppose that $G$ is a connected imperfect graph.  When $|V(G)|$ is small, the conclusion holds easily. Now we proceed with the inductive proof.

Without loss of generality, we may suppose that $G$ is a nonempty clique blowup of $F$. (Otherwise, $G$ is a nonempty clique blowup of some induced subgraph $F'$ of $F$, and we may write $F'$ as $F$.) Notice that $F$ is connected because the connectedness of $G$.

Let $u,v\in V(F)$. We use $Q_{v}$ to denote the clique of $G$ corresponding to $v$ of $F$, and let $Q_{u,v}=Q_u\cup Q_v$ if $uv\in E(F)$.

Let $T$ be a subset of $V(G)$ such that $|T\cap Q_v|=\min\{q, |Q_v|\}$ for any $v\in V(F)$. We define $$V_1=\bigcup_{uv\in E(F)}\{Q_{u,v}~:~\min\{|Q_v|, |Q_u|\}\leq q-1, \mbox{ and }\max\{|Q_v|, |Q_u|\}\geq \omega(G)-q+1\},$$
and let $V_2=V(G)\setminus V_1$.

Let $u_1v_1$ and $u_2v_2$ be two edges of $F$ such that $Q_{u_1,v_1}\cup Q_{u_2,v_2}\subseteq V_1$. Suppose that $|Q_{u_i}|\geq |Q_{v_i}|$ for $i\in \{1,2\}$. If $u_1u_2\in E(G)$, then $|Q_{u_1,u_2}|\geq 2\omega(G)-2q+2> \omega(G)$ since $\omega(G)\geq 2q+1$, a contradiction. Thus, $u_1u_2\notin E(F)$. Hence, $V_1\setminus T$ is a  union of pairwise anticomplete cliques of $G$. Let $V_1\setminus T=\cup_{i=1}^{\ell} K_i$, where $K_1,\cdots, K_{\ell}$ are such cliques. Clearly, $|K_i|\leq \omega(G)-1-q$ for $i\in\{1,\cdots, {\ell}\}$, and so we have that
\begin{equation}\label{w-1-q}
	\chi(G[V_1\setminus T])\leq\omega(G)-1-q.
\end{equation}

Let $uv$ be an edge of $F$ such that $Q_{u}\subseteq V_1$ and $Q_{v}\subseteq V_2$. If $|Q_u|\geq\omega(G)-q+1$, then $|Q_{v}|\leq\omega(G)-(\omega(G)-q+1)=q-1$, which implies that $Q_{v}\subseteq V_1$ and leads to a contradiction. Hence, $|Q_u|\leq q-1$, and so we can deduce that $V_1\setminus T$ is anticomplete to $V_2\setminus T$. Therefore,
\begin{equation}\label{chi-1}
	\chi(G-T)=\max\{\chi(G[V_1\setminus T]), \chi(G[V_2\setminus T])\}.
\end{equation}

\begin{claim}\label{clique number}
$\omega(G[V_2\setminus T])\leq \omega(G)-2q.$
\end{claim}
\pf For, $i\in \{1, 2\}$, let $F_i$ be a triangle-free graph such that $G[V_i]$ is a clique blowup of $F_i$. We will prove that for any maximal clique $K$ of $G[V_2]$, $|K\setminus T|\leq\omega(G)-2q$. Observe that either $K=Q_{x,y}$ for some $xy\in E(F_2)$, or $K=Q_z$ for some $z\in V(F_2)$.

Suppose that $K=Q_{x,y}$ for some $xy\in E(F_2)$. Since $Q_{x, y}\nsubseteq V_1$, we have that either $\min\{|Q_x|, |Q_y|\}\geq q$ or $\max\{|Q_x|, |Q_y|\}\leq \omega(G)-q$. If $\min\{|Q_x|, |Q_y|\}\geq q$, then $|K\setminus T|\leq \omega(G)-2q$. Thus we may assume that $\min\{|Q_x|, |Q_y|\}\leq q-1$ and $\max\{|Q_x|, |Q_y|\}\leq \omega(G)-q$. Without loss of generality, suppose that $|Q_x|\leq q-1$ and $|Q_x|\le |Q_y|\leq \omega(G)-q$. If $|Q_y|\geq q$, then $|K\setminus T|=|Q_y|-q\leq \omega(G)-2q$. If $|Q_y|\leq q-1$, then $|K\setminus T|=0<\omega(G)-2q$.

Next, we suppose that $K=Q_z$ for some $z\in V(F_2)$. Then, $z$ is an isolated vertex of $F_2$. Since $F$ is connected, $z$ must have a neighbor, say $z'$, in $V(F_1)$. If $|Q_{z'}|\geq\omega(G)-q+1$, then $|Q_z|\leq\omega(G)-|Q_{z'}|\leq q-1$, and so $Q_z\in V_1$, a contradiction. Hence, $|Q_{z'}|\leq q-1$. Since $z\in V(F_2)$, we have that $Q_z\nsubseteq V_1$, and so $|Q_z|\leq \omega(G)-q$. If $|Q_z|\geq q$, then $|K\setminus T|\leq \omega(G)-2q$. Otherwise, $|K\setminus T|=0<\omega(G)-2q$. This proves Claim~\ref{clique number}. \qed

Recall that $\chi(G')\leq\lceil\frac{p}{2q}\omega(G')\rceil$ if $G'$ is a clique blowup of $F$ with $\omega(G')\leq\max\{\frac{2q(p-q-2)}{p-2q}, 2q\}$. It follows that $\chi(G[T])\leq \lceil\frac{p}{2q}\omega(G[T])\rceil\leq p$ since $\omega(G[T])\leq 2q$. From $\omega(G)>\frac{2q(p-q-2)}{p-2q}$, we have that $\omega(G)-1-q\leq \lceil\frac{p}{2q}(\omega(G)-2q)\rceil$. By (\ref{w-1-q}), (\ref{chi-1}) and Claim~\ref{clique number}, and by the inductive hypothesis, we have that
\begin{eqnarray*}	
	\chi(G)&\leq&\chi(G-T)+\chi(G[T])\\
	       &\leq&\max\{\chi(G[V_1\setminus T]), \chi(G[V_2\setminus T])\}+p\\
	       &\leq&\max\{\omega(G)-1-q, \lceil\frac{p}{2q}(\omega(G)-2q)\rceil\}+p\\
	       &\leq&\lceil\frac{p}{2q}(\omega(G)-2q)\rceil+p\\
	       &=& \lceil\frac{p}{2q}\omega(G)\rceil.
\end{eqnarray*}
This completes the proof of Theorem~\ref{blowup}. \qed

\section{Proof of Theorem~\ref{main-1}}

In this section, we will prove Theorem~\ref{main-1}. Before that, we introduce a notation called {\em path extension} which was defined in \cite{WX2019} and plays an important role in our proof.

Let $n\geq3$ be an odd integer. For convenience, we use $[n]$ to denote the set $\{1, 2, \ldots, n\}$. Let $P=v_0v_1\cdots v_n$ be a path, and let $Q$ be obtained from blowing up each vertex $v_i$ into a nonempty clique $V_i$ for $0\leq i\leq n$. Let $r\geq \omega(Q)$ be an integer. Let $\phi$ be a coloring of $V_0\cup V_n$ with $r$ colors, and $S$ be a subset of $\phi(V_0)\setminus \phi(V_n)$ such that $|S|\leq\min\{|V_{2j}|~:~0\leq 2j\leq n\}$. It is easy to verify that, for each odd integer $j\in [n-1]$, we may greedily extend $\phi$ to an $r$-coloring $\phi'$ of $Q-V_j$ such that all colors of $S$ are used on vertices of $V_i$ while $i$ is even. We call $\phi'$ a ($Q, V_0, V_j, S$)-{\em path extension} of $\phi$, and call $V_0$ the {\em starting point} of the path extension.

The proof of Theorem~\ref{main-1} depends on the following Lemma~\ref{K_5}.

\begin{lemma}\label{K_5}
	Let $G$ be an (cap, even hole)-free graph. Then
\begin{itemize}
	\item [(i)] $\chi(G)\leq\lceil\frac{5}{4}\omega(G)\rceil$ if $G$ is $K_5$-free, and
	\item [(ii)]  $\chi(G)\leq\lceil\frac{7}{6}\omega(G)\rceil$ if $G$ is $(5$-hole, $K_{13})$-free.
\end{itemize}
\end{lemma}

\noindent\textbf{{\em Assuming Lemma~$\ref{K_5}$, we can prove Theorem}~\ref{main-1}}: Let $G$ be a (cap, even hole)-free graph. We prove that $\chi(G)\leq\lceil\frac{5}{4}\omega(G)\rceil$ by induction on $|V(G)|$. The result holds easily when $|V(G)|$ is small. Now, we proceed with the inductive proof and suppose $G$ is a connected imperfect graph. Moreover, we may assume that $G$ has no clique cutset and universal clique. Since $G$ has an odd hole, we may choose a maximal triangle-free induced subgraph $F$ of $G$ which has no clique cutset. By Lemma~\ref{cap}, we have that $G$ is a clique blowup of $F$.

Notice that $F$ is (even hole, triangle)-free. By Lemma~\ref{K_5}$(i)$, $\chi(G')\leq \lceil\frac{5}{4}\omega(G')\rceil$ for any clique blowup $G'$ of $F$ with $\omega(G')\leq4$. By applying Theorem~\ref{blowup} with $p=5$ and $q=2$, we have that $\chi(G)\leq \lceil\frac{5}{4}\omega(G)\rceil$. This proves the first conclusion of Theorem~\ref{main-1}.

If $G$ is further 5-hole-free, then  $\chi(G')\leq \lceil\frac{7}{6}\omega(G')\rceil$ for any clique blowup $G'$ of $F$ with $\omega(G')\leq12$ by Lemma~\ref{K_5}$(ii)$. By applying Theorem~\ref{blowup} with $p=7$ and $q=3$, we can prove the second conclusion of Theorem~\ref{main-1}. \qed

\medskip

Next, we will prove Lemma~\ref{K_5} in Subsection 4.1. From Lemma~\ref{odd hole} below, we see that both bounds of Lemma~\ref{K_5} are reachable.
\begin{lemma}\label{odd hole}{\em \cite{CX2024}}
Let $q\ge 2$, let $C=v_1v_2\cdots v_{2q+1}v_1$, and let $G$ be a clique blowup of $C$.
Then, $\chi(C^k)=\lceil\frac{2q+1}{2q}\omega(C^k)\rceil$ for every positive integer $k$,
and $\chi(G)\leq\lceil\frac{2q+1}{2q}\omega(G)\rceil$.
\end{lemma}

\subsection{Proof of Lemma~\ref{K_5}}

We first prove Lemma~\ref{K_5}$(i)$.
Suppose Lemma~\ref{K_5}$(i)$ does not hold. We may choose $G$ to be a minimum counterexample, that is $\chi(G)>\lceil\frac{5}{4}\omega(G)\rceil$, and any (even hole, cap, $K_5$)-free graph $H$ with order less than $G$ satisfies $\chi(H)\leq \lceil\frac{5}{4}\omega(H)\rceil$. By the minimality of $G$, we have that $G$ is a connected imperfect graph without clique cutset or universal clique, and $\delta(G)\geq \lceil\frac{5}{4}\omega(G)\rceil$.

Since $G$ has an odd hole, we may choose a maximal triangle-free induced subgraph graph $F$ of $G$ which has no clique cutset. Since $G$ is even hole-free, we have that $F$ is even hole-free of order at least 5 and is not a cube, and so $F$ is odd-signable. By Lemma~\ref{odd-signable}, $F$ can be obtained, starting from an odd hole, by a sequence of good ear additions. And $G$ can be constructed from $F$ by blowing up vertices of $F$ into nonempty cliques by Lemma~\ref{cap}.

If $F$ is an odd hole of length $2q+1$, then by Lemma~\ref{odd hole}, $\chi(G)\leq\max_{q\ge 2}\lceil\frac{2q+1}{2q}\omega(G)\rceil=\lceil\frac{5}{4}\omega(G)\rceil$, a contradiction. Therefore, $F$ is not an odd hole, and so we may assume that $F$ is obtained from $F'$ by adding a good ear, say $P_0=x_0x_1x_2\cdots x_m$, with attachments $x_0$ and $x_m$, from a hole $H$ of $F'$. Let $y$ be the
unique common neighbor of $x_0$ and $x_m$ in $F'$.  By the definition of good ear addition, $y$ has an odd number of neighbors in $P_0$, and so $y$ has a neighbor in $V(P_0)\setminus\{x_0, x_m\}$. We choose the first and second neighbor of $y$ in $P_0-x_0$, say $a$ and $b$ respectively. Since $F$ is (even
hole, triangle)-free, we have that $a$ must have an odd distance of at least 3 from $x_0$, and $b$ must have an odd distance of at least 3 from $a$. Hence, we may assume that $a=x_{2m_1-1}$ and $b=x_{2m_2}$ for some positive integers $m_1$ and $m_2$. Note that $2m_1-1\geq3$ and $2m_2\geq6$. Let $P=x_0x_1\cdots x_{2m_1-1}\cdots x_{2m_2}$, $P_1=x_0x_1\cdots x_{2m_1-1}$, and $P_2=x_{2m_1-1}x_{2m_1}\cdots x_{2m_2}$ be three segments of $P_0$.

\begin{figure}[htbp]
	\begin{center}
		\includegraphics[width=6cm]{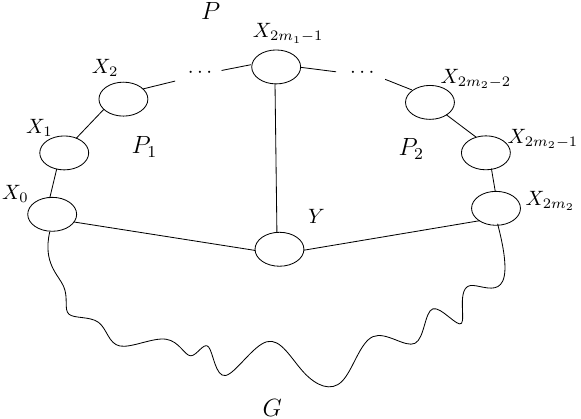}
	\end{center}
	\vskip -25pt
	\caption{Illustration of $G$, $P_1$ and $P_2$.}
	\label{fig-1}
\end{figure}

Notice that $G$ is obtained from $F$ by blowing up vertices into nonempty cliques. Let $Y$ be the clique of $G$ corresponding to $y$ of $F$, and
let $X_i$, for $0\leq i\leq 2m_2$, be the clique of $G$ corresponding to $x_i$ of $F$ (see Figure~\ref{fig-1}). Obviously, we have that
$$
\max\{|X_0|+|Y|, |X_{2m_1-1}|+|Y|, |X_{2m_2}|+|Y|, |X_{i}|+|X_{i+1}|, 0\leq i\leq 2m_2-1\}\leq \omega(G),
$$
and for each $i\in\{1,2,\cdots,2m_1-2, 2m_1,\cdots,2m_2-1\}$, and any $x_i\in X_i$,
$$
\lceil\frac{5}{4}\omega(G)\rceil\leq\delta(G)\leq d(x_i)=|X_{i-1}|+|X_i|+|X_{i+1}|-1\leq\min\{|X_{i-1}|+\omega(G)-1, |X_{i+1}|+\omega(G)-1\}.
$$

We can easily deduce that
$$
\lceil\frac{\omega(G)}{4}\rceil<|X_i|<\lfloor\frac{3\omega(G)}{4}\rfloor~\mbox{for}~0\leq i\leq2m_2,~\mbox{and}~|Y|<\lfloor\frac{3\omega(G)}{4}\rfloor.
$$
	
Since $2\leq\omega(G)\leq4$, we have that $|X_i|=2$ for $0\leq i\leq 2m_2$, and $1\leq|Y|\leq 2$. Therefore, $\omega(G)=4$, and $\chi(G)>5$. Since $G$ is a minimum counterexample, there exists a 5-coloring $\phi$: $V(G-\bigcup_{i=1}^{2m_2-1}X_i)\rightarrow [5]$, where $[5]=\{1,\cdots,5\}$.

Now we deduce a contradiction by using two steps of path extension to extend $\phi$ to $G$.
Without loss of generality, suppose $\phi(Y)\subseteq\{1, 2\}$ and $\phi(X_0)=\{3,4\}$. Extend $\phi$ to $G-\bigcup_{1\leq i\leq 2m_2-1, i\ne 2m_1-1}X_i$ as following. If $\phi(X_{2m_2})\ne\{4,5\}$, then we define $\phi(X_{2m_1-1})=\{4,5\}$; otherwise, we define $\phi(X_{2m_1-1})=\{3,5\}$. Now, it is easy to check that $|\phi(X_{2m_1-1})\setminus\phi(X_i)|\geq1$, where $i\in\{0,2m_2\}$. Choose $\alpha\in\phi(X_{2m_1-1})\setminus\phi(X_0)$ and $\beta\in \phi(X_{2m_1-1})\setminus\phi(X_{2m_2})$.
Let $Q_1=G[\bigcup_{i=0}^{2m_1-1}X_i]$ and $Q_2=G[\bigcup_{i=2m_1-1}^{2m_2}X_i]$. By making a ($Q_1, X_{2m_1-1}, X_2, \{\alpha\}$)-path extension of $\phi$, we get an extension $\phi_1$ of $\phi$ to $G-X_2-\bigcup_{2m_1\leq i\leq 2m_2-1}X_i$.

Since $|[5]\setminus (\phi(X_1)\cup \phi(X_3))|=5-|\phi(X_1)|-|\phi(X_3)|+|\phi(X_1)\cap \phi(X_3)|\geq 5-2-2+1= |X_2|$, we can color the vertices of $X_2$ and get an extension $\phi_2$ of $\phi_1$ to $G-\bigcup_{2m_1\leq i\leq 2m_2-1}X_i$. Similarly, by making a ($Q_2, X_{2m_1-1}, X_{2m_2-2}, \{\beta\}$)-path extension of $\phi_2$, we get an extension $\phi_3$ of $\phi_2$ to $G-X_{2m_2-2}$. Since $|[5]\setminus (\phi(X_{2m_2-1})\cup\phi(X_{2m_2-3}))|\geq |X_{2m_2-2}|$, we can color the vertices of $X_{2m_2-2}$ and finally extend $\phi_3$  to a 5-coloring of $G$, a contradiction. This completes the proof of Lemma~\ref{K_5}$(i)$.

\medskip

Now we turn to prove Lemma~\ref{K_5}$(ii)$. The approach is similar to the former case but somewhat complicated. Suppose Lemma~\ref{K_5}$(ii)$ does not hold. We still choose $G$  to be a minimum counterexample, i.e., $\chi(G)>\lceil\frac{7}{6}\omega(G)\rceil$, and any (even hole, cap, 5-hole, $K_{13}$)-free graph $H$ with order less than $G$ satisfies $\chi(H)\leq \lceil\frac{7}{6}\omega(H)\rceil$. Then, $G$ is connected, imperfect, has no clique cutset or universal clique, and $\delta(G)\geq \lceil\frac{7}{6}\omega(G)\rceil$. Let $\omega=\omega(G)$, and $k=\lceil\frac{7}{6}\omega\rceil$.

Moreover, with the same argument as in the proof of Lemma~\ref{K_5}$(i)$, let $F$ be a maximal triangle-free induced subgraph of $G$ which has no clique cutset and is not a cube. By Lemmas~\ref{odd-signable} and \ref{cap}, $F$ can be obtained from an odd hole of order at least 7 by a sequence of good ear additions, and $G$ is a nonempty clique blowup of $F$. It follows from Lemma~\ref{odd hole} that $F$ is certainly not an odd hole. Suppose $F$ is obtained from $F'$ by adding a good ear, say $P_0=x_0x_1\cdots x_m$, with attachments $x_0$ and $x_m$, and $y$ be the unique common neighbor of $x_0$ and $x_m$ in $F'$. Let $P$, $P_1$, $P_2$, $m_1$, $m_2$, $Y$, $X_i$, $0\leq i\leq 2m_2$, be defined as in the proof of above case $(i)$. Notice that $2m_1-1\geq5$ and $2m_2\geq10$ as $G$ is 5-hole-free. Also, we can easily deduce that, for each $i\in[2m_2-1]\setminus\{2m_1-1\}$ and any $x_i\in X_i$,
$$
\lceil\frac{7}{6}\omega\rceil\leq\delta(G)\leq d(x_i)=|X_{i-1}|+|X_i|+|X_{i+1}|-1\leq\min\{|X_{i-1}|+\omega-1, |X_{i+1}|+\omega-1\},
$$
and therefore
\begin{equation}\label{equ-1}
	\lceil\frac{\omega}{6}\rceil<|X_i|<\lfloor\frac{5\omega}{6}\rfloor~\mbox{for}~0\leq i\leq2m_2,~\mbox{and}~|Y|<\lfloor\frac{5\omega}{6}\rfloor.
\end{equation}

Since $\lceil\frac{5}{4}\omega\rceil=\lceil\frac{7}{6}\omega\rceil$ while $\omega\leq4$, we may assume that $5\leq\omega\leq12$ by Lemma~\ref{K_5}$(i)$. By the minimality of $G$, $G-\bigcup_{1\le i\le 2m_2-1} X_i$ admits a $k$-coloring, say $\phi$. In the following proofs, we will deduce a contradiction by extending $\phi$ to a $k$-coloring of $G$. The approach is similar to that of path extension. Roughly speaking, we choose a set $S$ of colors ({\em usually contained in or disjoint from $\phi(X_0)$}), and first extend $\phi$ to $\phi'$ by using the elements of $S$ to color some vertices of $X_i$ for each even $i$ (when $S\subseteq \phi(X_0)$), or to color some vertices of $X_i$ for each odd $i$ (when $S\cap \phi(X_0)=\emptyset$). Then, we further extend $\phi'$ to a $k$-coloring of $G$ by some greedy procedure.

\begin{claim}\label{1}
	For $0\leq i\leq 2m_2$, $|X_i|>\lceil\frac{\omega}{3}\rceil$.
\end{claim}
\pf Suppose it is not the case, and  assuming by symmetry that $|X_{i_0}|\leq\lceil\frac{\omega}{3}\rceil$ for some $i_0\in\{0, 1, \ldots, 2m_1-1\}$.
We first prove that
\begin{equation}\label{0,2m_1-1}
	i_0\notin\{0, 2m_1-1\}.
\end{equation}

\begin{figure}[htbp]
	\begin{center}
			\includegraphics[width=6cm]{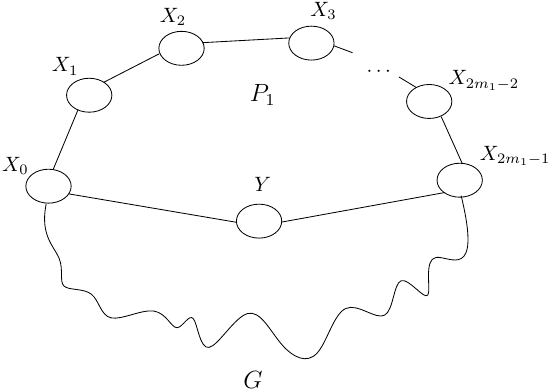}
		\end{center}
	\vskip -25pt
	\caption{Illustration of $G$ and $P_1$.}
	\label{fig-2}
\end{figure}

Suppose not, and assume without loss of generality that $i_0=0$. By our inductive hypothesis, let $\phi$ be a $k$-coloring of $G-\bigcup_{i=1}^{2m_1-2}X_i$ (see Figure~\ref{fig-2}). By (\ref{equ-1}), $|X_j|>\lceil\frac{\omega}{6}\rceil$, $0\leq j\leq 2m_1-1$, and so, we may choose $S\subseteq \phi(X_0)$ with $|S|=\lceil\frac{\omega}{6}\rceil$.

We first color $|S|$ vertices of $X_j$ with the colors in $S$ when $j$ is even in $\{2, 4, \ldots, 2m_1-4\}$, and greedily color all vertices of $X_{2m_1-2}$ using as many colors in $S$ as possible. Denote the resulted partial extension of $\phi$ as $\phi_1$. Then, we start from $X_{2m_1-4}$, and  extend $\phi_1$ to $\phi_2$  that colors sequentially all uncolored vertices of $X_j$ for even $j$ in $\{2m_1-4, \ldots, 2\}$, such that $|\phi_2(X_{j+2})\cap \phi_2(X_{j})|$ is maximum for each even integer $j$. Let $q$ be an odd integer in $\{3, \ldots, 2m_1-3\}$. Since $\max\{|X_{q-1}|, |X_{q+1}|\}\le \omega-|X_{q}|$, and since $|S|=\lceil\frac{\omega}{6}\rceil$, we have that $|\phi_2(X_{q-1})\cup\phi_2(X_{q+1})|\leq (\omega-|X_{q}|)+\lceil\frac{\omega}{6}\rceil$, and so
$|[k]\setminus (\phi_2(X_{q-1})\cup \phi_2(X_{q+1}))|\geq \lceil\frac{7\omega}{6}\rceil-(\omega-|X_{q}|)-\lceil\frac{\omega}{6}\rceil=|X_q|$. Since $|X_0|\leq\lceil\frac{\omega}{3}\rceil$, we have that $|[k]\setminus (\phi_2(X_{0})\cup \phi_2(X_2))|= \lceil\frac{7\omega}{6}\rceil-|\phi_2(X_0)|-|\phi_2(X_2)|+|\phi_2(X_0)\cap \phi_2(X_2)|\geq \lceil\frac{7\omega}{6}\rceil-\lceil\frac{\omega}{3}\rceil-(\omega-|X_{1}|)+
\lceil\frac{\omega}{6}\rceil \geq |X_{1}|$.

Hence, we can color all vertices of $X_q$ for each odd integer $q$ in $[2m_1-2]$, and further extend $\phi_2$ to a $k$-coloring of $G$. This completes the proof of (\ref{0,2m_1-1}).

\medskip

Now, we have $i_0\in [2m_1-2]$. Suppose that $i_0$ is odd. By (\ref{equ-1}), $|X_0|<\lfloor\frac{5\omega}{6}\rfloor$. Since $\lceil\frac{7\omega}{6}\rceil-(\lfloor\frac{5\omega}{6}\rfloor-1)>\lceil\frac{\omega}{6}\rceil$, we may choose $S\subseteq[k]\setminus\phi(X_0)$ such that $|S|=\lceil\frac{\omega}{6}\rceil$. Let $i_1=2$ if $i_0=1$, and let $i_1=i_0-1$ otherwise. Similar to the procedure of proving (\ref{0,2m_1-1}), we complete the proof of Claim~\ref{1} by extending $\phi$ to a $k$-coloring of $G$ as following.

Firstly, we color $|S|$ vertices of $X_j$ with the colors in $S$ when $j$ is odd in $[2m_1-3]$, and denote by $\phi_1$ the resulted partial extension of $\phi$.
Secondly, we start from $X_{2m_1-3}$, and  extend $\phi_1$ to $\phi_2$  that colors sequentially all uncolored vertices of $X_j$ for odd integer $j$ in $\{2m_1-3, \ldots, i_1\}$, such that $|\phi_2(X_{j+2})\cap \phi_2(X_{j})|$ is maximum for each odd integer $j$. Then, we start from $X_{1}$, and  extend $\phi_2$ to $\phi_3$  that colors sequentially all uncolored vertices of $X_h$ for odd integer $h$ in $[i_1]$, such that $|\phi_3(X_{h+2})\cap \phi_3(X_{h})|$ is maximum for each odd integer $h$.

Let $p$ be an even integer in $\{2, 3, \ldots, 2m_1-2\}$. If $p\ne i_1$, then since $\max\{|X_{p-1}|, |X_{p+1}|\}\le \omega-|X_{p}|$, we have that   $|\phi_3(X_{p-1})\cup\phi_3(X_{p+1})|\leq (\omega-|X_{p}|)+|S|=(\omega-|X_{p}|)+\lceil\frac{\omega}{6}\rceil$, and so $|[k]\setminus (\phi_3(X_{p-1})\cup \phi_3(X_{p+1}))|\geq \lceil\frac{7\omega}{6}\rceil-(\omega-|X_{p}|)-\lceil\frac{\omega}{6}\rceil=|X_p|$.
If $p=i_1$, note that either $i_1-1$ or $i_1+1$ equals $i_0$ and $|X_{i_0}|\le \lceil\frac{\omega}{3}\rceil$, then $|[k]\setminus (\phi_3(X_{p-1})\cup \phi_3(X_{p+1}))|\geq \lceil\frac{7\omega}{6}\rceil-(\omega-|X_{p}|)-\lceil\frac{\omega}{3}\rceil+
\lceil\frac{\omega}{6}\rceil\ge|X_p|$.
Therefore, one can easily extend $\phi_3$ to a $k$-coloring of $G$.

If $i_0$ is even, we may choose $S$ from $[k]\setminus\phi(X_{2m_1-1})$, and extend $\phi$ to a $k$-coloring of $G$ in almost the same way. This proves Claim~\ref{1}. \qed

\medskip
	
By Claim~\ref{1}, we have that
\begin{equation}\label{equ-2}
	\lceil\frac{\omega}{3}\rceil<|X_i|<\lfloor\frac{2\omega}{3}\rfloor~\mbox{for}~0\leq i\leq2m_2,~\mbox{and}~|Y|<\lfloor\frac{2\omega}{3}\rfloor.
\end{equation}
	
Recall that $5\leq\omega\leq 12$. From $|X_i|\geq \lceil\frac{\omega}{3}\rceil+1\geq 3$ for all $i$, we have that $6\leq\omega\leq12$. We will show, in the following three claims, that $\omega=12$.

\begin{claim}\label{2}
	$\omega\geq8$.
\end{claim}
\pf Suppose to its contrary that $\omega\leq7$. It follows from (\ref{equ-2}) that $\omega=6$, $|X_i|=3$ for all $i$, $1\leq|Y|\leq3$, and $k=\lceil\frac{7\omega}{6}\rceil=7$. Let $\phi$ be a $7$-coloring of $G-\bigcup_{i=1}^{2m_2-1}X_i$. Without loss of generality, we suppose that $\phi(Y)=[t]$ for some $t\in [3]$, and $\phi(X_0)=\{4,5,6\}$.

If $\phi(X_{2m_2})=\{4,5,6\}$, then let $\phi(X_{2m_1-1})=\{5,6,7\}$ and $S_1=\{1,7\}$; otherwise, let $\phi(X_{2m_1-1})=\{\alpha, \beta, 7\}$ for some $\alpha\in \phi(X_0)\setminus \phi(X_{2m_2})$ and $\beta\in [6]\setminus(\phi(Y)\cup \{\alpha\})$, and let $S_1=\{1, \alpha\}$. Now, $\phi$ is a 7-coloring of $G_1=G-\bigcup_{i\in[2m_2-1]\setminus\{2m_1-1\}} X_i$ such that
$$\{1, 7\}\cap\phi(X_{0})=\emptyset, \mbox{ and } S_1\cap\phi(X_{2m_2})=\emptyset.$$


Let $G_2=G-\bigcup_{1\le i\le m_1-1} X_{2i}-\bigcup_{2m_1\leq i\leq 2m_2-1} X_i$, and let $G_3=G-\bigcup_{1\le i\le m_2-1} X_{2i}$. We first extend $\phi$ to a 7-coloring $\phi_1$ of $G_2$ such that $\{1, 7\}\subseteq \phi_1(X_i)$ for odd $i$ in $[2m_1-3]$ and $\phi_1(X_{2m_1-3})\setminus\{1, 7\}\subseteq \phi(X_{2m_1-1})$. Then, we extend $\phi_1$ to a 7-coloring $\phi_2$ of $G_3$ such that $S_1\subseteq \phi_2(X_i)$ for odd $i$ in $\{2m_1+1, \ldots, 2m_2-1\}$ and $\phi_1(X_{2m_1+1})\setminus S_1\subseteq \phi(X_{2m_1-1})$.
Now, for each even integer $j$ in $\{2, \ldots, 2m_2-2\}$, $|\phi_2(X_{j-1})\cap \phi_2(X_{j+1})|\ge 2$, and so $|[7]\setminus \phi_2(X_{j-1})\cap \phi_2(X_{j+1})|=7-|\phi_2(X_{j-1})|-|\phi_2(X_{j+1})|+|\phi_2(X_{j-1})\cap \phi_2(X_{j+1})|\ge 3=|X_j|$. Therefore, we can extend $\phi_2$ to a 7-coloring of $G$, a contradiction. This proves Claim~\ref{2}. \qed

\medskip

Now, $8\leq\omega\leq12$, and $|X_i|\geq \lceil\frac{\omega}{3}\rceil+1\geq 4$ for all $i$ by (\ref{equ-2}). We prove next that $\omega\geq10$.

\begin{claim}\label{3}
	$\omega\geq10$.
\end{claim}
\pf Suppose to its contrary that $8\leq\omega\leq9$. By (\ref{equ-2}), we have that, $4\leq|X_i|\leq5$ for all $i$, $1\leq|Y|\leq5$, and $k=\lceil\frac{7\omega}{6}\rceil\in\{10, 11\}$. Let $\phi$ be a $k$-coloring of  $G-\bigcup_{i=1}^{2m_1-2}X_i$. We first show that
\begin{equation}\label{=5}
	\mbox{there exists an $i_0\in [2m_1-2]\setminus\{2,2m_1-3\}$ such that $|X_{i_0}|=5$.}
\end{equation}

Suppose not. Then, $|X_i|=4$ for each $i\in [2m_1-2]\setminus\{2,2m_1-3\}$. We choose $S\subseteq[k]\setminus\phi(X_0)$ with $|S|=2$, and color 2 vertices of  $X_j$ with the colors in $S$ for all odd $j\in [2m_1-3]$. Denote by $\phi'$ the resulted extension of $\phi$. Then, we start from $X_{2m_1-3}$, and extend $\phi'$ to $\phi''$ by sequentially coloring all uncolored vertices of $X_j$ for all odd $j\in [2m_1-3]$ such that $|\phi''(X_{j+2})\cap\phi''(X_j)|$ is maximum. Then, $|\phi''(X_{2m_1-3})\cap\phi''(X_{2m_1-1})|\ge |X_{2m_1-3}|-2$ which implies that
$|[k]\setminus (\phi''(X_{2m_1-3})\cup\phi''(X_{2m_1-1}))|\ge k-|X_{2m_1-3}|-(\omega-|X_{2m_1-2}|)+(|X_{2m_1-3}|-2)\ge|X_{2m_1-2}|$, and for all even $h$ in $[2m_1-3]$, $|\phi''(X_{h-1})\cap\phi''(h+1)|\ge 2$ which implies that
$|[k]\setminus (\phi''(X_{h-1})\cup\phi''(h+1))|\ge k-4-(\omega-|X_{h}|)+2=|X_{h}|$. Therefore, we can extend $\phi''$ to a $k$-coloring of $G$, a contradiction. This proves (\ref{=5}), and therefore, $\omega=9$ and $k=11$.

Now, we have $|X_{i_0-1}|=|X_{i_0+1}|=4$. From the symmetry between $X_0$ and $X_{2m_1-1}$,  we may suppose that $i_0$ is odd. Choose $S\subseteq\phi(X_0)$ with $|S|=4$, and choose $S_1\subseteq S$ with $|S_1|=2$. First, we extend $\phi$ to $\phi_1$ by coloring $X_j$ for even $j$ in $[i_0]$, such that $S\subseteq \phi_1(X_j)$, and coloring 2 vertices of $X_j$ for even $j$ in $\{i_0+1, \ldots, 2m_1-4\}$ with colors in $S_1$. Then, we extend $\phi_1$ to $\phi_2$ by first coloring $X_{2m_1-2}$ greedily, and sequentially coloring $X_j$ for even $j\in \{2m_1-4, \ldots, i_0+1\}$ such that $|\phi_2(X_j)\cap \phi_2(X_{j+2})|$ is maximum. Then, $|\phi_2(X_{h-1})\cap \phi_2(X_{h+1})|\ge 4$ when $h$ is odd in $[i_0-1]$, and $|\phi_2(X_{h-1})\cap \phi_2(X_{h+1})|\ge |X_{h-1}|-2$ when $h$ is odd in $\{i_0+1, \ldots, 2m_1-3\}$. It is certain that $|[k]\setminus (\phi_2(X_{i_0-1})\cup \phi_2(X_{i_0+1}))|= 11-|X_{i_0-1}|-|X_{i_0+1}|+|\phi_2(X_{i_0-1})\cap \phi_2(X_{i_0+1})|\ge 11-4-4+2=|X_{i_0}|$. For for each odd integer $h\in \{i_0+1, \ldots, 2m_1-3\}$, we have that
\begin{eqnarray*}	
	& &|[k]\setminus (\phi_2(X_{h-1})\cup \phi_2(X_{h+1}))|\\
	&=& 11-|\phi_2(X_{h-1})|-|\phi_2(X_{h+1})|+|\phi_2(X_{h-1})\cap \phi_2(X_{h+1})|\\
	&\geq& 11-|X_{h-1}|-(\omega-|X_{h}|)+(|X_{h-1}|-2)\geq |X_{h}|,
\end{eqnarray*}	
and for each odd integer $h\in [i_0-1]$, $|[k]\setminus (\phi_2(X_{h-1})\cup \phi_2(X_{h+1}))|\ge 11-|X_{h-1}|-|X_{h+1}|+|\phi_2(X_{h-1})\cap \phi_2(X_{h+1})|\ge 11-|X_{h-1}|-|X_{h+1}|+4\ge5\ge |X_{h}|$. Therefore, we can extend $\phi_2$ to a $k$-coloring of $G$, a contradiction. This proves Claim~\ref{3}. \qed

\medskip

By Claim~\ref{3} and (\ref{equ-2}), we have $10\leq\omega\leq12$, and $|X_i|\geq \lceil\frac{\omega}{3}\rceil+1\geq 5$ for all $i$. Now, we can prove that $\omega=12$. From now on, we let $G_1=G-X_{2}-\bigcup_{i=2m_1}^{2m_2-1}X_i$, let $G_2=G-\bigcup_{i=2m_1}^{2m_2-1}X_i$, and let $\phi$ denote a $k$-coloring of  $G-\bigcup_{i=1}^{2m_2-1}X_i$.

\begin{claim}\label{4}
$\omega=12$, and $\min\{|X_0|, |X_{2m_2}|\}=5$.
\end{claim}
\pf Suppose that Claim~\ref{4} does not hold, and thus either $10\leq \omega\leq11$ or $\min\{|X_0|,|X_{2m_2}|\}\geq6$. By (\ref{equ-2}), we have $5\leq |X_i|\leq7$ for all $i$, and $1\leq|Y|\leq7$.

If $10\leq \omega\leq11$, we have $|X_1|\leq 6\le \lceil\frac{\omega}{6}\rceil+4$ since $|X_2|\geq5$. If $\omega=12$, then  $|X_1|\le \omega-|X_0|\leq6=\lceil\frac{\omega}{6}\rceil+4$ since $|X_0|\geq 6$ by our assumption.
Therefore, we have by symmetry that
\begin{equation}\label{equ-3}
	|X_1|\le \lceil\frac{\omega}{6}\rceil+4,\mbox{ and } |X_{2m_2-1}|\le \lceil\frac{\omega}{6}\rceil+4.
\end{equation}

Since $|[k]\setminus(\phi(X_0)\cup \phi(Y))|\geq \lceil\frac{7\omega}{6}\rceil-\omega=\lceil\frac{\omega}{6}\rceil=2$, we can choose $T_1\subseteq[k]\setminus(\phi(X_0)\cup \phi(Y))$ with $|T_1|=2$. Since $|X_0|<\lfloor{2\omega\over 3}\rfloor$ by (\ref{equ-2}), we have $|[k]\setminus(T_1\cup \phi(X_0))|\geq \lceil\frac{7\omega}{6}\rceil-2-(\lfloor\frac{2\omega}{3}\rfloor-1)> 2$, and we can choose $T_1'\subseteq[k]\setminus(T_1\cup \phi(X_0))$ with  $|T_1'|=2$. Similarly, we may choose $T_2\subseteq[k]\setminus(\phi(X_{2m_2})\cup \phi(Y))$, and $T_2'\subseteq[k]\setminus(T_2\cup \phi(X_{2m_2}))$ such that $|T_2|=|T_2'|=2$. Let $S_1=T_1\cup T_1'$ and $S_2=T_2\cup T_2'$. So, $|S_1|=|S_2|=4$.

Let $\phi'$ be an extension of $\phi$ by coloring $X_{2m_1-1}$ such that $T_1\cup T_2\subseteq \phi'(X_{2m_1-1})$. Let $\phi_1$ be an extension of $\phi'$ to a $k$-coloring of $G_1$ in the following way: first color 4 vertices in $X_j$ for odd $j\in [2m_1-3]$ with colors in $S_1$, then starting from $X_{2m_1-2}$, sequentially color the uncolored vertices of $\bigcup_{i\in[2m_1-2]\setminus\{2\}} X_i$ greedily. This procedure is feasible since for each even $h\in[2m_1-2]\setminus\{2\}$, after coloring all vertices of $X_{h+1}$ and before $X_{h-1}$ has been fully colored,
\begin{eqnarray*}	
	& &|[k]\setminus (\phi'(X_{h-1})\cup \phi_1(X_{h+1}))|\\
	&=& \lceil\frac{7\omega}{6}\rceil-|\phi'(X_{h-1})|-|\phi_1(X_{h+1})|+|\phi'(X_{h-1})\cap \phi_1(X_{h+1})|\\
	&\geq& \lceil\frac{7\omega}{6}\rceil-4-(\omega-|X_{h}|)+2\geq |X_{h}|.
\end{eqnarray*}	
Now, let us consider $X_2$. Since $|X_1|\le \lceil\frac{\omega}{6}\rceil+4$ by (\ref{equ-3}), we have that
$|[k]\setminus (\phi_1(X_{1})\cup \phi_1(X_{3}))|=\lceil\frac{7\omega}{6}\rceil-|\phi_1(X_{1})|-|\phi_1(X_{3})|+|\phi_1(X_{1})\cap \phi_1(X_{3})|\ge\lceil\frac{7\omega}{6}\rceil-|X_1|-(\omega-|X_2|)+4\geq |X_{2}|$. So, we can extend $\phi_1$ to a $k$-coloring $\phi_2$ of $G_2$.

With almost the same argument as above, we may extend $\phi_2$  to $G$ by first coloring 4 vertices of $X_j$ with the colors in $S_2$ for odd $j\in \{2m_1+1, \ldots, 2m_2-1\}$, then sequentially and greedily coloring uncolored vertices of $X_{2m_1}$, $X_{2m_1+1}, \cdots, X_{2m_2-1}$. This contradicts the minimality of $G$, and completes the proof of Claim~\ref{4}. \qed

\medskip

Now, we have $\omega=12$, and we may assume by symmetry  that $|X_0|=5$. By (\ref{equ-2}), $5\leq |X_i|\leq 7$ for all $i$, $1\leq|Y|\leq7$, and $k=\lceil\frac{7\omega}{6}\rceil=14$. Recall that $\phi$ is a 14-coloring of $G-\bigcup_{i=1}^{2m_2-1}X_i$. Without loss of generality, we may suppose that $\phi(Y)=[t]$ for some $t\in [7]$, and $\phi(X_0)=\{8,9,10,11,12\}$. Next, we show that $|X_{2m_2}|=5$ also.

\begin{claim}\label{5}
	$|X_{2m_2}|=5$.
\end{claim}
\pf  Suppose to its contrary that $|X_{2m_2}|\geq6$. Then, $|\phi(Y)|\leq 6$, and so $7\notin\phi(Y)$. Let $T_1=\{7,13,14\}$ and $T_1'=\{1,2\}$. Since $|X_{2m_2}\cup Y|\le \omega=12$, we have that $|[14]\setminus(\phi(X_{2m_2})\cup \phi(Y))|\geq2$, and so can choose $T_2\subseteq [14]\setminus(\phi(X_{2m_2})\cup \phi(Y))$ with $|T_2|=2$. Let $T_2'\subseteq [14]\setminus(\phi(X_{2m_2})\cup T_2)$ and $|T_2'|=2$. Let $S_1=T_1\cup T_1'$ and $S_2=T_2\cup T_2'$. So, $|S_1|=5$ and $|S_2|=4$.

Since $|T_1\cup T_2|\leq5\leq |X_{2m_1-1}|$, we can extend $\phi$ to a 14-coloring $\phi'$ of $G-\bigcup_{i\in [2m_2-1]\setminus \{2m_1-1\}} X_i$ such that $T_1\cup T_2\subseteq \phi'(X_{2m_1-1})$.

Let $\phi_1$ be an extension of $\phi'$ to a $14$-coloring of $G_1$ in the following way: first color 5 vertices of $X_j$ for odd $j\in [2m_1-3]$ with colors in $S_1$, then starting from $X_{2m_1-2}$, sequentially color the uncolored vertices of $\bigcup_{i\in[2m_1-2]\setminus\{2\}} X_i$ greedily. This procedure is feasible since for each even $h\in[2m_1-2]\setminus\{2\}$,
$|[14]\setminus (\phi'(X_{h-1})\cup \phi_1(X_{h+1}))|=14-|\phi'(X_{h-1})|-|\phi_1(X_{h+1})|+|\phi'(X_{h-1})\cap \phi_1(X_{h+1})|\geq 14-5-(12-|X_{h}|)+3 = |X_{h}|$ before $X_{h-1}$ has been fully colored.
While for $X_2$, since $|X_2|\ge 5$, we have that
$|[14]\setminus (\phi_1(X_{1})\cup \phi_1(X_{3}))|=14-|\phi_1(X_{1})|-|\phi_1(X_{3})|+|\phi_1(X_{1})\cap \phi_1(X_{3})|\ge14-2(12-|X_2|)+5\geq |X_{2}|$. So, we can extend $\phi_1$ to a $k$-coloring $\phi_2$ of $G_2$.

Finally, we will extend $\phi_2$ to a 14-coloring $\phi_3$ of $G$ in the following way: first color 4 vertices of $X_j$ with the colors in $S_2$ for odd $j\in \{2m_1+1, \ldots, 2m_2-1\}$, then starting from $X_{2m_1}$, sequentially and greedily coloring uncolored vertices of $\bigcup_{i=2m_1}^{2m_2-1} X_{i}\setminus X_{2m_2-2}$ (the same reason as above, this can be done correctly). Since $|X_{2m_2}|\ge 6$, we have $|X_{2m_2-1}|\le 6$. Before finally coloring $X_{2m_2-2}$, we have $|[14]\setminus (\phi_3(X_{2m_2-1})\cup\phi_3(X_{2m_2-3}))|=14-|\phi_3(X_{2m_1-1})|-
|\phi_3(X_{2m_1-3})|+|\phi_3(X_{2m_1-1})\cap \phi_3(X_{2m_1-3})|\geq 14-6-(12-|X_{2m_2-2}|)+4=|X_{2m_2-2}|$. So, we can color $X_{2m_2-2}$ correctly, and extend $\phi_3$ to a 14-coloring of $G$. This contradicts the minimality of $G$, and completes the proof of Claim~\ref{5}. \qed

\medskip

By Claims~\ref{4} and \ref{5}, we have $\omega=12$, and $|X_{0}|=|X_{2m_2}|=5$. Choose $\alpha\in\phi(X_0)$ and $\beta\in \phi(X_{2m_2})$. Since $|[14]\setminus(\{\alpha, \beta\}\cup \phi(Y))|\geq 14-2-(12-|X_{2m_1-1}|)=|X_{2m_1-1}|$, we can extend the 14-coloring $\phi$ of $G-\bigcup_{i=1}^{2m_2-1}X_i$ to $G-\bigcup_{i\in [2m_2-1]\setminus\{2m_1-1\}} X_i$ by greedily coloring $X_{2m_1-1}$ with colors in $[14]\setminus(\{\alpha, \beta\}\cup \phi(Y))$.

Let $T_1\subseteq\phi(X_{0})\setminus \{\alpha\}$ and $T_2\subseteq\phi(X_{2m_2})\setminus \{\beta\}$ such that $|T_1|=|T_2|=2$, and let $S_1=\{\alpha\}\cup T_1$ and $S_2=\{\beta\}\cup T_2$. Clearly, $|S_1|=|S_2|=3$. We now extend $\phi$ to $\phi'$ by first coloring 3 vertices of $X_{i}$, for even $i\in [2m_1-4]$, with colors in $S_1$, and coloring all vertices of $X_{2m_1-2}$ such that $\{\alpha\}\in \phi'(X_{2m_1-2})$. Since for each odd integer $j\in \{2, \ldots, 2m_1-2\}$, $|[14]\setminus (\phi'(X_{j-1})\cup \phi'(X_{j+1}))| = 14-|\phi'(X_{j-1})|-|\phi'(X_{j+1})|+|\phi'(X_{j-1})\cap \phi(X_{j+1})|\geq 14-3-(12-|X_{j}|)+1=|X_{j}|$, we can extend $\phi'$ to a 14-coloring $\phi_1$ of $G-X_1-\bigcup_{2m_1\leq i\leq 2m_2-1} X_i$, by starting from $X_{2m_1-3}$ and sequentially coloring uncolored vertices of $\bigcup_{i=2}^{2m_1-2} X_i$. While for $X_1$, since $|X_0|=5$, we have $|[14]\setminus (\phi_1(X_{0})\cup \phi_1(X_{2}))|=14-|\phi_1(X_{0})|-|\phi_1(X_{2})|+|\phi_1(X_{0})\cap \phi_1(X_{2})|\geq 14-5-(12-|X_1|)+3=|X_{1}|$. Therefore, we can correctly coloring $X_1$ and further extend $\phi_1$ to a 14-coloring $\phi_2$ of $G-\bigcup_{2m_1\leq i\leq 2m_2-1}X_i$.

Take $X_{2m_2}, X_{2m_2-1}, \beta$ and $S_2$ acting the roles as $X_{0}, X_{1}, \alpha$ and $S_1$, respectively, we can extend $\phi_2$ to a 14-coloring of $G$ with the same argument as above that deals with $\bigcup_{i=0}^{2m_1-1} X_i$. This contradicts the minimality, and completes the proof of Lemma~\ref{K_5}$(ii)$ finally . \qed

\medskip

\noindent{\bf Remark:} Let $q\ge 3$ be an integer, and let ${\cal G}_q$ denote the set of (cap, even hole)-free graphs which do not contain an odd hole of length at most $2q-1$.  In this paper, we prove that every (cap, even hole)-free graph $G$ satisfies that $\chi(G)\leq\lceil\frac{5}{4}\omega(G)\rceil$, and every graph $G\in {\cal G}_3$ satisfies $\chi(G)\leq\lceil\frac{7}{6}\omega(G)\rceil$. These two bounds are reachable as shown by Lemma~\ref{odd hole}. It seems true that $\chi(G)\leq\lceil\frac{2q+1}{2q}\omega(G)\rceil$ for every graph $G\in {\cal G}_q$ for any integer $q\geq3$. By Theorem~\ref{blowup}, to prove this statement, it suffices to verify that $\chi(G)\leq\lceil\frac{2q+1}{2q}\omega(G)\rceil$ for graphs $G\in {\cal G}_q$ with $\omega(G)\le \max\{\frac{2q(p-q-2)}{p-2q}, 2q\}$.
While using our argument of proving Lemma~\ref{K_5}, it may contain too many detailed local analysis and the proof would be too boring. It would be nice if someone can prove the statement, or present some counterexample.

\end{document}